\documentclass[12pt]{article}
\usepackage{natbib}
\usepackage{hyperref}
\usepackage{amssymb}
\usepackage{dsfont}
\usepackage{hhline}
\usepackage{epsfig}
\begin{document}
\newtheorem{theorem}{Theorem}
\newtheorem{remarque}{Remark}
\newtheorem{lemma}{Lemma}
\newtheorem{corollaire}{Corollary}

%
%
\vskip 3mm

\vskip 3mm \baselineskip=15pt

\noindent SOME UNIFORM LIMIT RESULTS IN  ADDITIVE REGRESSION MODEL

\vskip 5mm \noindent Mohammed DEBBARH  \vskip 2mm \noindent

\noindent Universit\'e Paris 6

\noindent 175,\ Rue du Chevaleret,  75013 Paris.

\noindent debbarh@ccr.jussieu.fr.
 \vskip 4mm \baselineskip=18pt

\noindent ABSTRACT

We establish some uniform limit results in the setting of additive
regression model estimation. Our results allow to give an
asymptotic 100\% confidence bands for these components. These
results are stated in the framework of i.i.d random vectors when
the marginal integration estimation method is used.
\vskip 4mm

\section{Introduction}
For $ d \geq 2 $, let $({\bf X}, Y)$ be an $\mathbb{R}^d \times
\mathbb{R}$-valued random vector. The regression function of $Y$
given ${\bf X}= {\bf x} $ is defined, for any $ {\bf x} \in
\mathbb{R}^d $, by
\begin{equation}\label{regression}
m({\bf x}) = E( Y | {\bf X} = {\bf x}).
\end{equation}
For $ 0 < p < \infty $, the optimal $L_{p}$ rate of convergence of
a nonparametric estimate of $m$ is of order $n^{\frac {-k}{2k+d}}$
when $m$ is assumed to be a $k$-times differentiable function and
for $p = \infty $, the optimal rate is $(n^{-1} \log n)^{\frac
{k}{2k+d}}$ (see, \cite{Stone}). This rate of convergence which
depends on the dimension $d$ of the covariable ${\bf X}$ becomes
worse as the dimensionality of the problem increases. In the
literature, this phenomena is known under the name of ``curse of
dimensionality''. To reduce the dimension impact upon the
estimates, Stone (1985) proposed several sub-models  of model
$(\ref{regression})$. More particularly, he studied the
nonparametric additive regression model in which the multivariate
regression function is written as the sum of univariate functions,
i.e,
\begin{equation}\label{additive}
 m({\bf x}):= m_{add}({\bf x}) =  \mu + \sum_{l=1}^{d} m_{l}(x_{l}).
\end{equation}
To study the model $(\ref{additive})$, several estimation methods
have been proposed, in the literature. We cite, the method based
on B-spline (see, \cite{Stone}), the method based on the
backfitting algorithm (see, \cite{Hastie}); hereafter, we make use
of the marginal integration method, (see, \cite{Newey},
\cite{Auestad}and Linton and Nielsen (1995)). The additive
regression components have motivated the work of many researchers,
we refer to \cite{Camlong-Viot-Sarda-and-Vieu2000} for a survey on
the asymptotic normality  of the additive components under a
mixing condition and to \cite{Sperlich2003} for nonparametric
estimation and testing of integration in additive models.
\cite{debbarh2004} established the law of iterated logarithm for
additive regression model components under the independence
assumption on the sample $({\bf X}_i, Y_i)_{i=1,...,n}$.
 In this paper, we establish some
 uniform limit results in probability in an additive regression model.
 Similarly to results stated by
\cite{Deheuvels2004} for functionals of a distribution based on
kernel-type estimation, our results allow us to build asymptotic
100\%
 confidence bands for the components we estimate.
\section{Main results}
First,  introduce some  notations and assumptions. We denote by
$f$ and $f_l $ the joint density of ${\bf X}$ and
 the marginal density of $X_l$, for $l=1,...,d$, respectively. We consider the following
 assumptions upon $m$, $f$ and $f_l, l=1,...,d$. The functions $f$ and $f_l$ are  continuous
with compact supports and there exist  $ \ b,\ b_{l},\ B,\
\mbox{and} \ B_l$ such that\\
\\
$(F.1)$ \ $0 < b \leq f({\b x}) \leq B < \infty \ \mbox{and} \ 0 <
b_l \leq f_l (x_l) \leq B_l < \infty.$ \\
$(F.2)$ \ $m$ \mbox{ is
$k$-times continuously differentiable}.\\
$(F.3)$ \ $f$ \mbox{ is
$k'$-times continuously differentiable}, $k'>kd$.\\
\\
Throughout, $(\ell_n)_{n \geq 1}$, $(h_n )_{n\geq 1} \ \mbox{and}
\ (h_{l,n})_{n\geq 1},\ 1 \leq l \leq d,$ are
sequences of positive constants satisfying the following conditions\\
\\
$(H.1)$  \ $h_{l,n} \rightarrow 0, \mbox{and} \ \ell_n \rightarrow
0 \ \mbox{as} \ n \rightarrow
\infty$.\\
$(H.2)$ \ $nh_{1,n} \rightarrow +\infty$, $n h_{1,n} /\log(n)
\rightarrow \infty \ \mbox{as} \ n
\rightarrow \infty$.\\
$(H.3)$ \ $ n h_{1,n} h_{1,n}^{2k_1}...h_{d,n}^{2k_d}/|\log h_{1,n}|  \longrightarrow 0, \ k_1+...+k_d=k, \ \mbox{and} \ h_{1,n} \log(n)/\ell_{n}^d |\log h_{1,n}| \longrightarrow 0 \ \mbox{as} \ n \rightarrow \infty. $ \\
$(H.4)$ \ $h_n \sim n^{-2k/(2k+1)}$ and $nh_{1,n}h_n^2/|\log h_{1,n}| \rightarrow 0$ as $n \rightarrow \infty$.\\
\\
Let $I= \prod_{i=1}^d I_i = \prod_{i=1}^d [a_i, c_i]~\mbox{and}~
J= \prod_{i=1}^d J_i = \prod_{i=1}^d [a_i', c_i']$ be two fixed
pavements of $\mathbb{R}^d$ such that $a_i'< a_i<c_i<c_i',
~1\leq i \leq d.$ Furthermore, we consider the following assumptions upon the random variable $Y$.\\
\\
$(M.1)$ \ $Y \mathbb{I}_{\{{\bf X} \in J\}}$ is bounded.\\
\\
Let $({\bf X}_i, Y_i )_{i=1,...,n}$ be a $n$-sample with the same
 distribution as $({\bf X}, Y)$. Let $L$ be a kernel on $\mathbb{R}^d$,  of order $k'$, bounded and with compact
support. We define the kernel estimator $\widehat f_n $ of the
density $f$ by
\begin{eqnarray*}
{\widehat f}_n ({\bf x}) = \frac{1}{n \ell_{n}^d} \sum_{i=1}^{d}
L\Big(\frac{{\bf x} - {\bf X}_i }{\ell_n } \Big).
\end{eqnarray*}
To estimate the multivariate regression function defined in
(\ref{regression}), we will be used  the tow following estimators,
\begin{eqnarray}\label{estimateur}
\widehat m_{n}({\bf x}) = \sum_{i=1}^{n} {\frac{Y_{i}}{n {\widehat
f}_n (X_{i})}} \Big( \prod_{l=1}^{d} \frac{1}{h_{l,n}}
K_{l}\Big({\frac{x_{l} - X_{i,l}}{h_{l,n}}}\Big) \Big),
\end{eqnarray}
and
\begin{eqnarray}\label{estimateur2}
\widetilde m_{n}({\bf x}) = \sum_{i=1}^{n} {\frac{Y_{i}}{nh_n^d
{\widehat f}_n (X_{i})}} K\Big({\frac{{\bf x} - {\bf
X}_{i}}{h_{n}}}\Big),
\end{eqnarray}
 where, the kernel functions $K$ and $
K_{l},\ l=1,...,d,$ are bounded, continuous, 
and integrate to one. In addition, we assume that $K_l$ satisfies the following conditions\\
\\
$(K.1)$\ $\mathbb{K}:=\prod_{l=1}^d K_l$ is of order $k$.\\
\\
$(K.2)$\ $K_l(u) = 0$ for $u \not\in [-\frac{\lambda_l}{2},
\frac{\lambda_l}{2}),$ for some $0< \lambda_l < \infty.$\\
 \\
$(K.3)$\ $\mathbb{K}=\prod_{l=1}^d K_l$ is square integrable
function in the linear span (the set of finite linear
combinations) of functions $\Psi \geq 0$ satisfying the following
property: the subgraph of $\Psi$, $\{(s,u):\Psi(s)\geq u\}$, can
represented  as a finite number of Boolean operations among sets
of the form $\{(s,u):p(s,u)\geq \phi(u)\}$, where $p$ is a
polynomial on $\mathbb{R}^d\times \mathbb{R}$ and $\phi$ arbitrary
real function.\\ In particular this is satisfied by $K({\bf
x})=\phi(p({\bf x}))$, $p$ being a polynomial and $\phi$
a bounded real function of bounded variation (see \cite{Nolan} and \cite{gine2002}).\\
\\
As already mentioned, the marginal integration method  will be
used to estimate the additive components (see \cite{Linton1995}
and \cite{Newey}). Towards this aim, for all ${\bf
x}=(x_1,...,x_{d})\in\mathbb{R}^d$ and every  ${\bf
x}_{-l}=(x_1,..,$ $x_{l-1},x_{l+1},..,x_d)$,  $l=1,...,d$, set
$q({\bf x} ) = \prod_{l=1}^{d} q_l (x_l)$ and $q_{-l}({\bf x}_{-l}
) = \prod_{j\neq l} q_j (x_j)$. Then, the $l$-th component $\eta_l
$ of the additive model is given by
\begin{eqnarray}
\eta_{l}(x_{l}) = \int_{\mathbb{R}^{d-1}} m({\bf x}) q_{-l}({\bf
x}_{-l}) d{\bf x}_{-l} - \int_{\mathbb{R}^d} m({\bf z}) q({\bf z})
d{\bf z},\quad l=1,...,d, \label{additive_component}
\end{eqnarray}
in such a way that the two following equalities hold,
\begin{eqnarray}
&&\eta_{l}(x_l) = m_{l}(x_l) - \int_{\mathbb{R}} m_{l}(z)
q_l(z)dz, \quad l=1,...,d, \label{relation_additive_component}\\
&&m({\bf x})= \sum_{l=1}^d \eta_{l}(x_l) + \int_{\mathbb{R}^d}
m({\bf z})q({\bf z}) d{\bf z}
\label{additive_component_marginale}.
\end{eqnarray}
In view of (\ref{relation_additive_component}) and
(\ref{additive_component_marginale}), $\eta_{l}$ and $m_{l}$ are
equal up to an additive constant. Therefore, $\eta_{l}$ is an
additive component too, fulfilling a different identifiability
condition. Note also that $\eta_{l}=m_{l}$ for the particular
choice $q_l=f_l$, $l=1,...,d$. However, $f_l$ is generally
unknown, and $\eta_{l}\ne m_{l}$ in practice. From
(\ref{estimateur}) and (\ref{additive_component}), a natural
estimate of the $l$-th component $\eta_{l}$ is given by
\begin{eqnarray}\label{enq1}
\widehat \eta_l(x_{l}) = \int_{\mathbb{R}^{d-1}}\!
\widehat{m}_{n}({\bf x}) q_{-l}({\bf x}_{-l}) d{\bf x}_{-l} -
\int_{\mathbb{R}^d}\! \widehat{m}_{n}({\bf z}) q({\bf z}) d{\bf
z},\ l=1,...,d.
\end{eqnarray}
From (\ref{additive_component_marginale}) and (\ref{enq1}), we
derive an estimate $\widehat m_{add}$ of the additive regression
function,
\begin{eqnarray}
\widehat m_{add} ({\bf x})& = & \sum_{l=1}^d \widehat
\eta_{l}(x_l) + \int_{\mathbb{R}^d} \widetilde{m}_{n}({\bf z})
q({\bf z})d{\bf z} \label{estim_add}.
\end{eqnarray}
In other respects, we impose the following assumptions on the
known  integration density functions $q_{-l}$ and $q_l$,
$l=1, ..., d$,\\
\\
$(Q.1)~~ q_{-l} \mbox{ is bounded and
continuous, } l=1, ..., d.$\\
$(Q.2)~~ q_l \mbox{ has } k \mbox{ continuous and bounded
derivatives, with compact support } \mathcal{C}_l \subset I_l,
~l=1, ..., d.$\\
\\
Let $\phi(u_1)$ be a continuous function on the interval $I_1$
defined by
$$\phi(u_1 ) = \int_{\mathbb{R}^{d-1}} \frac{H(\bf{u})}{f({\bf u}_{-1}|u_1 )}
q_{-1} ({\bf u}_{-1}) d{\bf u}_{-1}, $$
 where $$ H({\bf u}) = E(Y^2| {\bf X} = {\bf u}), \ \ {\bf u}=(u_1,...,u_d)\in \mathbb{R}^d.$$
 Consider the following quantity
\begin{eqnarray}\label{def_var}
\sigma_{l} = \sup_{x_{l}\in I_l} \sqrt{
\frac{\phi(x_{l})}{f_{l}(x_{l}) } \int_{\mathbb{R}} K_{l}^2}.
\end{eqnarray}
The following results describe the asymptotic behavior  of the
estimates $\widehat \eta_1$ and $\widehat m_{add}$. From now on, $
\stackrel{P}{\rightarrow}$  denote the convergence in probability.
\begin{theorem}\label{th1}
Under the hypotheses $(F.1)-(F.3)$, $(H.1)-(H.4)$, $(K.1)-(K.3)$,
$(Q.1)-(Q.2)$ and
 $(M.1)$, we have, as $n \longrightarrow \infty$,
\begin{eqnarray}
\sqrt{\frac{nh_{1,n}}{2|\log h_{1,n} |}}\sup_{x_1 \in I_1} \pm \{
\widehat \eta_1 (x_{1}) -  \eta_1 (x_1 ) \}
\stackrel{P}{\longrightarrow} \sigma_1. \label{eq_theo1}
\end{eqnarray}
\end{theorem}
Theorem \ref{th2} below is valid under the additional condition
that for any all $ 1 \leq l \leq d$, \ $h_{l,n} = h_{1,n}$.
\begin{theorem}\label{th2}
Under the hypotheses $(F.1)-(F.3)$, $(H.1)-(H.4)$, $(K.1)-(K.3)$,
 $(Q.1)-(Q.2)$ and
 $(M.1)$, we have, as $n \longrightarrow \infty$,
$$\sqrt{\frac{nh_{1,n}}{2|\log h_{1,n}|}}\sup_{{\bf x}\in I} \pm \Big\{\widehat{m}_{add}({\bf x})
 - m({\bf x})\Big\}\stackrel{P}{\longrightarrow} \sum_{l}^{d}\sigma_l.
$$
\end{theorem}
\section{Application}
\subsection{Confidence bands}
Let $\sigma_{1,n}(x_1 )$ be the  estimator of $\sigma_1 (x_1 )$,
with $\sigma_1 (x_1 ) = \sqrt{\phi(x_1 )/f_1 (x_1 )}$. We will
consider the data-dependent function $L_{n} (x_1 )$, defined by,
\begin{eqnarray*}
L_{n} (x_1 ) = \Big\{ \frac{2 \log(1/h_{1,n})}{nh_{1,n} }\times
\sigma_{1,n} (x_1 )\Big\}^{1/2} \Big[ \int_{\mathbb{R}} K_{1}^{2}
\Big]^{1/2},
\end{eqnarray*}
where
\begin{eqnarray*}
\sigma_{1,n}^2 (x_1)= \frac{1}{nh_{1,n}}\sum_{i=1}^n Y_i^2
K_1\Big(\frac{x_1-X_{i,1}}{h_{1,n}}\Big) \int_{\mathbb{R}^{d-1}}
\frac{ \prod_{l \neq 1} \frac{1}{h_{l,n}} K_{l}\Big({\frac{x_{l} -
X_{i,l}}{h_{l,n}}}\Big) }{\widehat f_n^2({\bf x})} q_{-1}({\bf
x}_{-1})~d{\bf x}_{-1}.
\end{eqnarray*}
We obtain asymptotic simultaneous confidence bands for
$\eta_{1}(x_1)$  in the following sense. For each $0 < \epsilon <
1$, we have
\begin{eqnarray*}
\mathbb{P} \Big\{ \eta_{1}(x_{1}) \in {\Big[ \widehat \eta_{1}(x_1
) - (1+\epsilon) L_{n}(x_1 ), ~~ \widehat \eta_{1} (x_1 ) +
(1+\epsilon) L_{n} (x_1 ) \Big]}, \ \forall x_1  \in I_1 \Big\}
\longrightarrow 1,
\end{eqnarray*}
and
\begin{eqnarray*}
\mathbb{P} \Big\{ \eta_1 (x_1 ) \in \Big[ \widehat \eta_1 (x_1 ) -
(1- \epsilon ) L_n (x_1 ),~~ \widehat \eta_{1} (x_1 ) + (1 -
\epsilon ) L_n (x_1 ) \Big], \ \forall x_1 \in I_1 \Big\}
\longrightarrow 0.
\end{eqnarray*}
We say then that the intervals
\begin{eqnarray*}
 \Big[ A_{n,1} (x_1 ),~~ B_{n,1} ( x_1 ) \Big] = \Big[ \widehat \eta_{1}
(x_1 ) - L_n (x_1 ), ~~ \widehat \eta_{1} (x_1 ) + L_n ( x_1 )
\Big],
\end{eqnarray*}
 provide asymptotic simultaneous confidence bands (at an
asymptotic confidence level of 100 \%) for $\eta_1(x_1) \
\mbox{over} \  x_1 \in I_1 .$ We deduce the asymptotic confidence
bands for $m_{add}$, over ${\bf x} \in I,$
\begin{eqnarray*}
\Big[A_n({\bf x}),~~B_n({\bf x})\Big]= \Big[\mu_n+\sum_{l=1}^d
A_{n,l}(x_l),~~ \mu_n+\sum_{l=1}^d
B_{n,l}(x_l)\Big],~\mbox{where}~ \mu_n=\int_{\mathbb{R}^d}
\widehat m_n({\bf x}) q({\bf x})d{\bf x}.
\end{eqnarray*}
Following  a general statistical practice, for finite values of
the sample size $n$, we recommend the use of the asymptotic 100 \%
confidence bands.  Our results do not provide  confidence regions
in the usual sense, since they are not based on a specified level
$1- \alpha $. Instead, they hold with probability tending to $ 1 \
\mbox{as} \ n \rightarrow \infty.$
\section{Proofs}
Let $\mathcal{G}$ be a class of  pointwise measurable functions
satisfying conditions ($\mathcal{C}$) in the Appendix. We denote
by $\alpha_{n}(.)$ the multivariate empirical processus based upon
$ ({\bf X}_{1}, Y_{1}), ({\bf X}_{2}, Y_{2}),...$ \ and indexed by
the class of functions $ \mathcal{G}$. More precisely, $\alpha_n $
is defined for $g \in \mathcal{G}$ by,
$$ \alpha_{n} (g) = \frac{1}{\sqrt{n}} \sum_{i=1}^{n}
\Big(g({\bf X}_{i}, Y_{i}) - E( g ( {\bf X}_{i}, Y_{i}) ) \Big).
$$ For any real valued function $ \phi $ defined on a set $B$, we
use the notation $||\phi ||_{B} = \sup_{x  \in B} | \phi(x) |:= ||
\phi ||$. Recalling that  $I_1 = [a_1, c_1]$, let $0 < \eta < 1 $
be  a fixed  number and set, for $n \geq 1,$
\begin{equation}\label{discritisation}
x_{1,j} = a_1 + j \eta h_{1,n}, \ \ 0 \leq j \leq l_n :=
\Big[\frac{c_1 -a_1}{\eta h_{1,n}} \Big],
\end{equation}
where $[u]$ denotes the integer part of $u$. For ${\bf
X}_i=(X_{i,1},...,X_{i,d}),$  $1 \leq i \leq n$, set
\begin{equation}\label{var}
 g_{n}^{x_1} ({\bf X}_i, Y_i) = Y_i G({\bf X}_i) K_{1} \Big( \frac{x_{1}
- X_{i,1}}{h_{1,n}}\Big),
\end{equation}
where,
\begin{eqnarray}\label{G}
 G({\bf X}_i) = \frac{1}{ f({\bf X}_i)}
\int_{\mathbb{R}^{d-1}}  \Big( \prod_{l \neq 1} \frac{1}{h_{1,n} }
K_{l} \Big(\frac{x_{l} - X_{i,l}}{h_{l,n}}\Big) q_{l} (x_{l})\Big)
d{\bf x}_{-l}.
\end{eqnarray}
For $ n \geq 1 $ and any $ 0 \leq j \leq l_n $, let $
\mathcal{G}_{n} = \{ g_{n}^{x_{1,j}}: \ 0 \leq j \leq l_n \} $.
Obviously, for each $0 \leq j \leq l_n$ and any $x_1 \in I_1,$ we
have
$$ || g_{n}^{x_{1,j}} || + || g_{n}^{x_1} || \leq \kappa, \ \ \mbox{where}\ \ \kappa \ \mbox{is a positive constant}.$$
In the first time, we  assume that the density $f$  of  ${\bf X}$
is known. Let ${\widehat {\widehat m}_n}$ be the estimator of the
regression function when $f$ is known,
$$\widehat {\widehat m}_{n}({\bf x}) = \frac{1}{n} \sum_{i=1}^{n} \frac {Y_{i}}{f({\bf X}_{i})} \Big( \prod_{l=1}^{d} \frac{1}{h_{l,n}} K_{l}\Big(\frac {x_{l} - X_{i,l}}{h_{l,n}}\Big)\Big). $$
Using the estimator $(\ref{enq1})$ of the additive regression
model components, we obtain
\begin{eqnarray*}
nh_{1,n} \Big( \widehat {\widehat \eta}_{1}(x_{1}) - E(\widehat
{\widehat \eta}_{1}(x_{1})) \Big)
& = & nh_{1,n} \int_{\mathbb{R}^{d-1}} \Big(\widehat {\widehat m}_{n}({\bf x}) - E(\widehat {\widehat m}_{n} ({\bf x}))\Big) q_{-1} ({\bf x}_{-1}) d{\bf x}_{-1} \\
                                                               &   & - nh_{1,n} \int_{\mathbb{R}^{d}} \Big(\widehat {\widehat m}_{n}({\bf x}) - E(\widehat {\widehat m}_{n} ({\bf x}))\Big) q (x) d{\bf
                                                               x}, \\
& = & \sqrt{n} \alpha_{n} (g_{n}^{x_1}) - \int_{\mathbb{R}}
\sqrt{n} \alpha_{n} (g_{n}^{x_1}) q_1 (x_1) dx_1,
\end{eqnarray*}
where,
\begin{eqnarray}\label{alphag}
\sqrt{n} \alpha_{n} (g_{n}^{x_1}) = nh_{1,n}
\int_{\mathbb{R}^{d-1}} \Big(\widehat {\widehat m}_{n}({\bf x}) -
E(\widehat {\widehat m}_{n} ({\bf x}))\Big) q_{-1} ({\bf x}_{-1})
d{\bf x}_{-1}.
\end{eqnarray}
The proof of Theorem \ref{th1} is based on a number of additional
lemmas.
\begin{lemma} \label{lem1} Assume that the conditions $(F.1)-(F.3)$, $(H.1)-(H.3)$, $(K.1)-(K.3)$,
 $(Q.1)-(Q.2)$ and
 $(M.1)$, are satisfied. Then,  we have, as $n \longrightarrow \infty$
$$ \frac {\sup_{x_1 \in I_1} \pm \alpha_{n} (g_{n}^{x_{1
}})  }{\sqrt {2 h_{1,n} |\log  h_{1,n}|}} \stackrel{P}{
\longrightarrow}  \sigma_{1}.$$
\end{lemma}
Following \cite{Deheuvels2004} and \cite{Mason2000}, the proof of
Lemma \ref{lem1} is split up into two part. First, we establish
the upper bound, afterward, we state the lower bound.
\subsection{Proof of Lemma \ref{lem1}}
\subsubsection{Upper bound part} The main tools used in the proof
are the discretization and the properties of empirical
process oscillations.\\
\\
{\bf Part~1.}~~ We examine the behavior of our process on an
appropriately chosen grid, $(x_{1,j})_{1 \leq j \leq l_n}$ of
$I_1$. Assume that the assumptions of Lemma \ref{lem1} hold. It
follows from (\ref{var}) that
\begin{eqnarray}\label{ingalit}
Var\Big( g_{n}^{x_{1,j}} ({\bf X}_i, Y_i)) \Big) & \leq & E\Big( \Big(g_{n}^{x_{1,j}} ({\bf X}_i, Y_i)\Big)^2 \Big) \nonumber\\
                                     & \leq & E\Big( G^{2} ({\bf X}_i ) Y_i^{2} K_1^2 \Big( \frac{x_{1,j} - X_{i,1}}{h_{1,n}} \Big)\Big) \nonumber \\
                                     & \leq & E\Big( E( Y_{i}^{2} | {\bf X}_i ) G^{2} ({\bf X}_i ) K_1^2 \Big( \frac{x_{1,j} - X_{i,1}}{h_{1,n}} \Big)\Big) \nonumber \\
                                     & \leq & E\Big( H({\bf X}_i ) G^{2} ({\bf X}_i ) K_1^2 \Big( \frac{x_{1,j} - X_{i,1}}{h_{1,n}} \Big)\Big).
\end{eqnarray}
But, making use of classical changes of variables and Taylor
expansion, we get under $(K.1)$ with ${\bf
h}_{-1}=(h_{2,n},...,h_{d,n})^T$ and $0 < \theta <1$,
\begin{eqnarray}
&&\int_{\mathbb{R}^{d-1}} \prod_{l \neq 1 } \Big(\frac{1}{h_{l,n}}
K_{l} \Big(\frac {x_{l} - u_{l}}{h_{l,n}} \Big) q_{l} (x_{l})\Big)
d{\bf x}_{-1} \nonumber \\& = &
\int_{\mathbb{R}^{d-1}}  \prod_{l \neq 1 }  \Big( K_{l} (v_{l}) q_{l} (v_{l} h_{l,n} + u_{l}) \Big)dx_{-1} \nonumber\\
              & = & \int_{\mathbb{R}^{d-1}}  \Big( \prod_{l \neq 1 }   K_{l} (v_{l}) \Big) q_{-1}({\bf v}_{-1} {\bf h}_{-1}+ {\bf u}_{-1}) d{\bf v}_{-1} \nonumber\\
              & = & \int_{\mathbb{R}^{d-1}} \!\! \Big( \prod_{l \neq 1 }   K_{l} (v_{l}) \Big) \Big[q_{-1} ({\bf u}_{-1})
              + \!\!\!\!\!\!\!\! \sum_{k_2+...+k_d=k} v_{2}^{k_2}... v_{d}^{k_d}
h_{2,k}^{k_2}...h_{d,k}^{k_d} \frac{\partial^k q_{-1}}{\partial
v_{2}^{k_1}...\partial v_{d}^{k_d} } ({\bf v}_{-1} {\bf h}_{-1}
\theta +
{\bf u}_{-1}) \Big] d{\bf v}_{-1} \nonumber \\
              & = & q_{-1} ({\bf u}_{-1}) + o(1). \nonumber
\end{eqnarray}
Therefore, it hold that,
\begin{equation}\label{equat}
\Big(\int_{\mathbb{R}^{d-1}} \prod_{l \neq 1 } \Big(
\frac{1}{h_{1,n}} K_{l} \Big(\frac {x_{l} - u_{l}}{h_{l,n}} \Big)
q_{l} (x_{l})\Big) d{\bf x}_{-1} \Big)^2  =   q_{-1}^{2} ({\bf
u}_{-1}) +  o(1).
\end{equation}
Moreover, in view of the (\ref{G}), we have
\begin{eqnarray}\label{esp}
&  & E\Big( H({\bf X}_i ) G^{2} ({\bf X}_i ) K_1^2 \Big(
\frac{x_{1,j} - X_{i,1}}{h_{1,n} } \Big)\Big) \nonumber \\& =  &
\int_{\mathbb{R}^d} \frac{H({\bf u})}{f({\bf u})} \Big(
\int_{\mathbb{R}^{d-1}} \prod_{l \neq 1 }  \frac{1}{h_{1,n}} K_{l}
\Big(\frac {x_{l} - u_{l}}{h_{l,n}} \Big) q_{l} (x_{l})
dx_{-1}\Big)^2 K_{1}^{2}\Big( \frac{x_{1,j} - u_{1}
}{h_{1,n}}\Big) d{\bf u},
\end{eqnarray}
Combining (\ref{ingalit}), (\ref{equat}) and (\ref{esp}), we
obtain
\begin{eqnarray}
Var\Big(g_{n}^{x_{1,j}} ({\bf X}_i, {\bf Z}_i, \delta_i) \Big)
\leq \sigma_{1}^2 h_{1,n} + o(h_{1,n}).
\end{eqnarray}
Applying Bernstein's inequality to the sequence of random
variables,
$$ Z_{r} = g_{n}^{x_{1,j}} ({\bf X}_{r}, Y_{r}) - E(g_{n}^{x_{1,j}} ({\bf X}_{r}, Y_{r})), \ \ r = 1,...,n
,$$ we obtain, for $n$ large enough, that,
\begin{eqnarray}
&  & \mathbb{P} \Big\{ \max_{ 1 \leq j \leq  l}  | \alpha_{n} (
g_{n}^{x_{1,j} })|  \geq \sigma_1 (1 + \rho ) \sqrt{2
h_{1,n}|\log h_{1,n}| } \Big\} \nonumber\\
& \leq & 2 (  l + 1 )  \exp\left( \frac{-2 \sigma_{1}^{2} (1 + \rho )^2 h_{1,n} |\log h_{1,n}| }{ 2 \sigma_{1}^2 h_{1,n} + \frac{2 M}{3 \sqrt{n}} \sigma_{1} \sqrt{2 h_{1,n} |\log h_{1,n}| }} \right) \nonumber\\
& \leq & 2( l + 1) h_{1,n}^{1 + \frac{\rho}{2}}.
\end{eqnarray}
\\
{\bf Part~2.}~~ Under assumption $(K.3)$, the class of functions
$$
\Big\{ K\Big(\frac{{\bf t} - .}{a}\Big): {t} \in \mathbb{R}^d, a
\in  \mathbb{R}^d \setminus\{0\}  \Big\}
$$
is a bounded VC class of measurable functions. Now, consider the
class
$$\mathcal{F}=\Big\{ b y   K\Big(\frac{{\bf t} - .}{a}\Big): {t} \in \mathbb{R}^d, a
\in  \mathbb{R}^d \setminus\{0\}, |b|\leq D \Big\},$$ where $D>0$
is the bound of the function $y G({\bf x})$. Arguing exactly as in
pages 254 and 255 of Deheuvels and Mason (2004), one can show that
$\mathcal{F}$ fulfills $\mathcal{C}$. An easy argument now shows
that
$$\mathcal{G}'=\Big\{ b y   K\Big(\frac{{\bf t} - .}{a}\Big)-  b' y   K\Big(\frac{{\bf t'} - .}{a'}\Big): {\bf t}, {\bf t}' \in \mathbb{R}^d, a, a' \in \mathbb{R} \setminus\{0\}, |b|\leq D, |b'|\leq D
\Big\}$$ fulfills conditions $\mathcal{C}$. As over $\mathbb{G}$
function we can take
\begin{eqnarray*}
\mathbb{G}({\bf u},v)=2Cv||K||_{\infty}.
\end{eqnarray*}
Since $\mathcal{G}_n' \subset \mathcal{G}'$. We study the behavior
of our process between the grid points $x_{1,j},~ x_{1,j+1},$ with
$1 \leq j \leq l_n.$ Toward this aim, consider for $0 \leq j \leq
l_n,$ the class of functions
$$ \mathcal{G}_{n,j}^{'} = \{
g_{n}^{x_{1,j}} - g_{n}^{x_1 },\ x_{1,j} \leq x_1 \leq x_{1,j+1}
\}~~\mbox{and}~~ \mathcal{G}_{n}^{'} = \bigcup_{j}
\mathcal{G}_{n,j}^{'}.$$ There exists an absolute constant $B$,
such that for any $\epsilon
> 0 $, one can find a $ \eta_{\epsilon}$ such that whenever $(\ref{discritisation})$
holds, with $ 0 < \eta < \eta_{\epsilon},$ we have,
\begin{eqnarray}\label{R0}
\mathbb{P} \Big\{||n^{1/2} \alpha_{n} ||_{\mathcal{G}_{n}^{'}}
\geq B \sqrt{\epsilon n h_{1,n} |\log h_{1,n}|} \Big\} = o(1).
\end{eqnarray}
Indeed,  we see that, uniformly over $g \in \mathcal{G}_{n,j}^{'}
\subset \mathcal{G}_{n}^{'},$\ $|| g || \leq \kappa.$ Moreover, by
similar arguments as those used  in the proof of
$(\ref{ingalit})$, we have,
\begin{eqnarray}\label{R1}
\sigma_{\mathcal{G}_{n}^{'}}^{2} = \sup_{g \in
\mathcal{G}_{n}^{'}} {\rm Var}(g ({\bf X}, Y)) \leq 4h_{1,n}
\sigma_{1}^2.
\end{eqnarray}
Therefore by Fact 1 (see the Appendix), for all $t > 0$ we have
for suitable finite constants   $A_{1}, A_{2} > 0$,
\begin{eqnarray}\label{R3}
&& \mathbb{P} \Big\{ || n^{1/2 } \alpha_{n}
||_{\mathcal{G}_{n}^{'}} \geq A_{1} \Big( E|| \sum_{i=1}^{n}
\epsilon_{i} g({\bf X}_{i}, Y_{i}) ||_{\mathcal{G}_{n}^{'}} + t
\Big) \Big\} \nonumber \\& \leq & 2 \Big\{ \exp\Big(- \frac{A_{2}
t^2}{n \sigma_{\mathcal{G}_{n}^{'}}^2}\Big) + \exp \Big( - \frac
{A_{2}
t}{\kappa}\Big) \Big\} \nonumber \\
& \leq & 2 \Big\{ \exp\Big(- \frac{A_{2} t^2}{2n h_{1,n}
\sigma_{1}^2 }\Big) + \exp \Big( - \frac {A_{2} t}{\kappa}\Big)
\Big\}.
\end{eqnarray}
Next, by using $(\ref{discritisation})$, in combination with Fact
2 (see the Appendix), we obtain that
\begin{equation}\label{R4}
E||\sum_{i=1}^{n} \epsilon_{i} g_{n}^{x_1} ({\bf X}_{i}, Y_{i})
||_{\mathcal{G}_{n}^{'}} \leq A^{'} \sqrt{\nu n h_{1,n}  \log
(1/h_{1,n})}.
\end{equation}
Where $A^{'}$ is an absolute constant. Thus,  using $(\ref{R4})$,
we get from $(\ref{R3})$ that
\begin{eqnarray}\label{RRR}
& & \mathbb{P} \Big \{||n^{1/2} \alpha_{n}
||_{\mathcal{G}_{n}^{'}}
\geq 2 A_{1} A^{'} \sqrt{\nu n h_{1,n} |\log h_{1,n}|} \Big\}\\
& \leq & 2  \exp \Big( - \frac{A_{2} (A_{1} A^{'} )^2 \nu |\log
h_{1,n}|}{\sigma_{1}^2} \Big) + 2 \exp\Big( - \frac{A_{2} (A_{1}
A^{'}) \sqrt{\nu n h_{1,n} | \log h_{1,n}| }}{M} \Big) \nonumber
\\ & = & o(1).
\end{eqnarray}
Taking $ B = 2 A_{1} A^{'} \sqrt{\nu \epsilon} $ in $(\ref{RRR})$,
we  complete the proof of (\ref{R0}).\\
\\
For $1 \leq j \leq l_n,$ $\mathcal{G}_{n,j}^{'} \subseteq
\mathcal{G}_{n},$  we see that
\begin{equation}\label{ineq2}
\frac{\max_{0 \leq j \leq l} ||n^{1/2} \alpha_n
||_{\mathcal{G}_{n,j}^{'} }}{\sqrt{2nh_{1,n}|\log h_{1,n}|}} \leq
 \frac{ ||n^{1/2} \alpha_n ||_{\mathcal{G}_{n}^{'}}}{\sqrt{2nh_{1,n}|\log h_{1,n}|}}.
\end{equation}
Using the statement (\ref{R0}) and the inequality (\ref{ineq2})
with $A = \frac{B}{\sigma_1 \sqrt{2}}$, we obtain
\begin{equation}\label{oscillation}
\mathbb{P}\Big\{\frac{\max_{1 \leq j \leq l_n} ||n^{1/2}
\alpha_n||_{ \mathcal{G}_{n,j}^{'}}}{\sqrt{2n h_{1,n} |\log
h_{1,n}|}} \geq \sigma_1 A \sqrt{\epsilon} \Big\}=o(1),
\end{equation}
\\
{\bf Conclusion }: By combining  (\ref{discritisation}) and
(\ref{oscillation}), we conclude that there  exists an absolute
constant $A>0$, such that
\begin{eqnarray*}
\mathbb{P} \Big\{  \frac {\sup_{{x_{1}} \in I_1}|
\alpha_{n}(g_{n}^{x_{1} }) | } {\sqrt{2  h_{1,n} |\log h_{1,n}| }}
> ( 1 + \rho + A \sqrt{ \epsilon } ) \sigma_{1} \Big\}
& \leq &  \mathbb{P} \Big\{ \max_{1 \leq j \leq l }  \frac { |\alpha_{n} ( g_{n}^{x_{1,j}}) | }{\sqrt{2h_{1,n} |\log h_{1,n}| }} > ( 1  + \rho ) \sigma_{1} \Big\} \\
&     & +  \mathbb{P} \Big\{  \max_{1 \leq j \leq l_n } \frac {||
n^{\frac{1}{2}} \alpha_{n} ||_{\mathcal{G}_{n,j}'}}{\sqrt{2n
h_{1,n} |\log h_{1,n}| }} > A \sqrt{\epsilon } \sigma_{1} \Big\}.
\end{eqnarray*}
Since for any $ \epsilon > 0$, we can choose $\rho > 0$ and $
\epsilon > 0 $ small enough so that $\rho + A \sqrt{\epsilon} <
\epsilon$. We obtain the upper bound result in the case where $f$
is known,
$$ \mathbb{P} \Big\{ \frac {\sup_{x_1 \in I_1}|\alpha_{n}(g_{n}^{x_{1}})|}{\sqrt{2  h_{1,n} |\log h_{1,n}| }} > ( 1 + \epsilon) \sigma_1 \Big\} = o(1). $$
\subsubsection{Lower bound part }
In order to prove lower bound result, we gather hereafter some
technical results (see, Einmahl and Mason (2000)). Let $Z_{1},
Z_{2},.....$ be a sequence of i.i.d random vectors
 taking values in $\mathbb{R}^d$. For each $n \geq 1$, consider the
empirical distribution function defined by,
$$ G_{n} ({\bf s}) = n^{-1} \sum_{i=1}^{n} \mathbb{I}_{ \{ {\bf Z}_{i} \leq {\bf s} \} }, \ s \in \mathbb{R}^{d+1} ,$$
where as usual ${\bf z} \leq {\bf s}$ means that each component of
${\bf z}$ is less  than or equal to the corresponding component of
$s$. For any measurable real valued function $g(.)$ defined on
$\mathbb{R}^{d+1}$, set
$$ G_{n} (g)= \int_{\mathbb{R}^{d}} g({\bf s})\   dG_{n}({\bf s}), \ \ \ \mu(g) = E(g({\bf Z}))  \ \mbox{and}  \  \sigma (g) = \sqrt{Var ( g({\bf Z}))}.$$
Let $\{ a_{n}: n \geq 1\}$ denote a sequence of positive constants
converging to zero and  satisfying the condition $|\log(a_n )|/
\log\log(n) \rightarrow \infty$. For some sequence of integer
number $k_n$, consider a sequence of sets of real valued
measurable functions on $\mathbb{R}^{d+1}$, $\mathcal{G}_{n} = \{
g_{i}^{(n)};\ i = 1,...,k_{n} \}$, defined by the following
conditions:
\\
$ \noindent {\bf (a)} ~~~~ \mathbb{P}  \Big( g_{i}^{(n)}(Z) \neq 0
, \  g_{i}^{(n)} (Z) \neq 0 \Big) = 0, \ \ \forall 1 \leq i \neq j
\leq k_{n}\ \mbox{and} \sum_{i=1}^{k_{n}} P\{g_{i}^{(n)}(Z) \neq
0) \leq 1/2.
$\\
Furthermore, assume that for some $0 < r <\infty$,\\
$ {\noindent ({\bf b})} ~~~~ a_{n} k_{n} \rightarrow r \ \
\mbox{as} \ \ n \rightarrow \infty.
$\\
For some $ - \infty < \mu_{1} < \mu_{2} < \infty $ and  $0 <
\sigma_{1} < \sigma_{2} < \infty,$ uniformly in $i = 1,...,k_{n},$
we have for  $n$ large enough,\\
$ {\noindent ({\bf c})} ~~~~ a_{n} \mu_{1} \leq \mu ( g_{i}^{(n)}
) \leq a_{n} \mu_{2}\  \ \mbox{and}\ \ \sqrt{ a_{n}} \sigma_{1}
\leq \sigma ( g_{i}^{(n)} ) \leq \sqrt{a_{n}} \sigma_{2}.
$\\
 For some $0 < M_1 <
\infty$, uniformly in $i= 1,...,k_{n}$, we have for  $n$ large
enough, \\
$
{\noindent ({\bf d})} ~~~~|| g_{i}^{(n)} || \leq M_1.$\\
\\
The following lemma due to Einmahl and Mason (2000) is the main
tool to prove our result. We will work only in the "+" case, the
arguments for the "-" case can be obtained similarly.
\begin{lemma}\label{lem2}
Under assumptions $({\bf a})-({\bf d})$, for each $0 < \epsilon <
1$, we have
$$ \mathbb{P} \Big\{ \max_{1 \leq i \leq k_n } \frac{n^{1/2}
\{G_n(g_{i}^{(n)}) - \mu(g_{i}^{(n)})
\}}{\sigma(g_{i}^{(n)})\sqrt{2 |\log(a_n )|}} \geq 1 - \epsilon
\Big\} \rightarrow  1.$$
\end{lemma}
\noindent{{\it Proof:}}   \ See Proposition 2 of Einmahl and Mason (2000).\\
\\
In order to apply the result of Lemma \ref{lem2}, we need to check
the validity of the conditions $({\bf a})-({\bf d})$ in our
setting. For any $ \epsilon
> 0$, select a sub-interval ${\it J_1 = [a_1', c_1']} \ \mbox{of} \ I_1 =
[a_1, c_1]$, such that
$$ \inf_{u_1 \in J_1} \sqrt{\frac{\phi(u_1 )}{f_1 (u_1)}} \Big[\int_{\mathbb{R}}
K_1 ^2 \Big]^{1/2} > \sigma_1 (1 - \epsilon/2) $$ and
$$ {\bf P} \{ X_1 \in J_1\} \leq 1/2.$$
Consider the following points in the interval  $J_1$
$$ x_{1,j} = {\it a_1'} + 2jh_{1,n}, \ \mbox{for} \ j = 1,...,
[(b_1'-a_1')/2h_{1,n}]-1:= k_n .$$ Then, it is easy to see that
the condition $({\bf b})$ is satisfied with $a_n = h_{1,n}$, i.e.
$$\lim_{n \rightarrow \infty} h_{1,n}k_n \approx [(b_1 -a_1)/2].$$ For
each $x_{1,j}$, $ 1\leq j \leq k_n , $ define the functions
$$ g_{j}^{(n)} ({\bf x}, y) = y G ({\bf x} ) K_{1} \Big( \frac{x_{1,j} - x_{1}}{h_{1,n}}
\Big),$$ it follows from $(F.1),$ that
$$ || g_{j}^{(n)}|| \leq M ||G||\times  ||K||:= M_1,$$
so the condition $({\bf d})$ is verified.\\
\\
Now, we verify the validity of the condition $({\bf a})$.  To this
end, recall that $K_1$ is compactly supported. therefore,
\begin{eqnarray*}
 g_{i}^{(n)} ({\bf X}, Y) & \neq & 0 \Longleftrightarrow \Big|\frac{x_{1,i} - X_{i,1}}{h_{1,n}}
 \Big|\leq \frac{1}{2},
\end{eqnarray*}
and then
\begin{eqnarray*}
 |x_{1,j} - X_{i,1}| & = & |x_{1,j} - x_{1,i}+ x_{1,i} - X_{i,1}| \geq
 2h_{1,n} -  \frac{h_{1,n}}{2}.
\end{eqnarray*}
Hence ,
 for $1 \leq i \neq j \leq k_n,$
$$\mathbb{P} \Big\{\  g_{i}^{(n)} ({\bf X}, Y) \neq 0 \
\mbox{and}\  g_{j}^{(n)} ({\bf X}, Y) \neq 0 \Big\} = 0.$$ To
check the validity of  the condition $({\bf c})$, recall the
inequality $(\ref{ingalit})$ and observe that, for $x_1 \in J_1$,
\begin{eqnarray*}
{\rm Var}( g_{i}^{(n) } ({\bf X}, Y)) & =    & E\Big( Y_{i}^{2} G^{2} (X_{i} ) K_{1}^{2} \Big( \frac{x_{1} - X_{i,1} }{h_{1,n}} \Big)\Big) -  E\Big( Y_{i} G (X_{i} ) K \Big( \frac{x_{1} - X_{i,1} }{h_{1,n}} \Big)\Big)^2 \\
& =    & \int_{\mathbb{R}} K_{1}^{2} \Big( \frac{x_{1} - u_{1}}{h_{1,n}}\Big) \frac{\phi(u_{1})}{f_{1} (u_{1})} du_{1} + o(h_{1,n}) \\
& \geq & h_{1,n} \Big(\int K_{1}^{2}\Big) \inf_{u_1 \in J_1} \frac{\phi(u_{1})}{f_{1} (u_{1})}  + o(h_{1,n}) \\
& \geq & \sigma_{1}^{2} (1 - \epsilon/2) h_{1,n}.
\end{eqnarray*}
Now, we can use the result of Lemma \ref{lem2}, which yields,
\begin{eqnarray*}
&   & \mathbb{P} \Big\{ \max_{1 \leq i \leq k_{n}} \frac{n^{1/2}
\{ G_{n}( g_{i}^{(n)} ) - \nu(g_{i}^{(n)}) \}}{\sqrt{2 {\rm Var}
(g_{i}^{(n)} )
|\log h_{1,n}|}} < 1 - \epsilon \Big\} \\
&\leq & \mathbb{P} \Big\{\max_{1 \leq i \leq k_{n}} \frac{n^{1/2}
(G_{n}(g_{i}^{n}) - \nu(g_{i}^{(n)} ))}{\sqrt{2 h_{1,n} |\log
h_{1,n} |}} < (1 - \epsilon)^{3/2} \sigma_1 h_{1,n}^{1/2} \Big\}.
\end{eqnarray*}
For $a_{n} = h_{1,n}$, we get from last inequality
$$\mathbb{P} \Big\{ \max_{1 \leq i \leq k_{n}} \frac{n^{1/2} (G_{n}(g_{i}^{n}) - \nu(g_{i}^{(n)} ))}{\sqrt{2 h_{1,n} |\log h_{1,n} |}} < (1 - \epsilon)^{3/2} \sigma_1 h_{1,n}^{1/2} \Big\} = o(1).$$
By setting $(1 - \epsilon )^{3/2} = 1 - \frac{\epsilon'}{2}$, and
using the inequality,
$$\frac{\sup_{x_1 \in I_1} \sqrt{n}\alpha_{n} (g_{n,1}^{x_{1}}) }{\sqrt{2 n h_{1,n} |\log h_{1,n}|}} \geq
\max_{1 \leq i \leq k_{n}} \frac{n^{1/2} (G_{n}(g_{i}^{n}) -
\nu(g_{i}^{(n)} ))}{\sqrt{2 h_{1,n} |\log h_{1,n} |}},$$ we obtain
$$  \mathbb{P} \Big\{ \frac{\sup_{x_1 \in I_1} \alpha_n( g_{n}^{x_1 }) }{\sqrt{2h_{1,n},n |\log(h_{1,n} )|}} < (1- \frac{\epsilon}{2}) \sigma_1
\Big\}= o(1).$$ Combining the results of Part 1 and Part 2, it
follows that,
$$ \Big| \frac{\sup_{x_1 \in I_1}  \alpha_{n} (g_{n}^{x_{1}})}{\sqrt{2  h_{1,n} |\log h_{1,n}|}} - \sigma_1 \Big| = o_{\bf P}(1). $$
Similarly one may show that
$$\Big| \frac{\sup_{x_1 \in I_1} - \alpha_{n} (g_{n}^{x_{1}})}{\sqrt{2  h_{1,n} |\log h_{1,n}|}} - \sigma_1 \Big| = o_{\bf P}(1). $$
This completes the proof of Lemma \ref{lem1}.$\square$

\subsection{Proof of Theorem \ref{th1}}
Let us use the decomposition,
\begin{eqnarray*}
\pm\{\widehat \eta_1(x_1) - \eta_1(x_1)\}&=&\pm\{\widehat
\eta_1(x_1)- \widehat {\widehat \eta}_1(x_1)\}\pm
\{\widehat {\widehat \eta}_1(x_1)- E(\widehat {\widehat \eta}_1(x_1))\}\pm\{E(\widehat {\widehat \eta}_1(x_1))-\eta_1(x_1)\}\\
&=& T_1(x_1)+T_1(x_2)+T_3(x_3).
\end{eqnarray*}
First, consider the term $T_1(x_1)$. It hold that,
\begin{eqnarray*}
\widehat m_{n}({\bf x})  = \widehat {\widehat m}_{n}({\bf x}) +
\frac{1}{n} \sum_{i=1}^{n} Y_{i} \Big(
\prod_{i=1}^{d}\frac{1}{h_{l,n}} K_{l}\Big(\frac {x_{l} -
X_{i,l}}{h_{l,n}}\Big)\Big ) \frac{f({\bf X}_{i}) - \widehat f_{n}
({\bf X}_{i})}{f({\bf X}_{i}) \widehat f_{n}({\bf X}_{i})}.
\end{eqnarray*}
It follows that,
\begin{eqnarray}
\Big|\widehat m_{n}({\bf x})  -  \widehat {\widehat m}_{n}({\bf
x}) \Big|
& \leq & \frac{1}{n} \sum_{i=1}^{n} |Y_{i}| \Big| \prod_{i=1}^{d}
\frac{1}{h_{l,n}} K_{l}\Big(\frac {x_{l} - X_{i,l}}{h_{l,n}}\Big)
\Big| \times \frac{\sup_{1\leq i\leq n} \Big|f({\bf X}_{i}) -
\widehat f_{n} ({\bf X}_{i})\Big|}{\Big|f({\bf X}_{i})\widehat
f_{n}({\bf X}_{i}) \Big|}.
\end{eqnarray}
Using, for example, the following result due to  \cite{Ango-Nze},
\begin{eqnarray}\label{angonze}
\sup_{{\bf x}\in I} \mid \widehat f_{n}({\bf x}) - f({\bf x}) \mid
= O\Big( \sqrt{\frac {\log (n)}{n\ell_{n}^d}}\Big)~~~\mbox{a.s.},
\end{eqnarray}
we obtain under the assumptions $(F.1)$, $(M.1)$ and
$(K.1)-(K.2)$, for $n$ large enough,
\begin{eqnarray}
\sup_{x\in  I} \mid \widehat m_{n}({\bf x}) - \widehat {\widehat
m}_{n}({\bf x}) \mid = O\Big( \sqrt{\frac
{\log(n)}{n\ell_n^d}}\Big) ~~~\mbox{a.s.}.
\end{eqnarray}
Observe that,
\begin{eqnarray*}
&& \hspace{-1.5cm}\Big|\widehat {\widehat  \eta}_{1}(x_ {1}) -
\widehat \eta_{1}(x_ {1}) \Big|\\& = & \Big|
\int_{\mathbb{R}^{d-1}} \Big[ \widehat {\widehat m}_{n}({\bf x}) -
\widehat m_{n}({\bf x}) \Big] q_{-1}({\bf x}) d{\bf x}_{-1} +
\int_{\mathbb{R}^d}
\Big[ \widehat {\widehat m}_{n}({\bf x}) - \widehat m_{n}({\bf x}) \Big] q({\bf x}) d{\bf x} \Big|, \\
&\leq & \int_{\mathbb{R}^{d-1}} \sup_{{\bf x}\in I} \Big| \widehat
{\widehat m}_{n}({\bf x}) - \widehat m_{n}({\bf x}) \Big|
q_{-1}({\bf x}_{-1})
dx_{-1} +  \int_{\mathbb{R}^d} \sup_{{\bf x}\in I} \Big| {\widehat {\widehat m}}_{n}({\bf x}) - \widehat m_{n}({\bf x})\Big| q({\bf x}) d{\bf x}, \\
& \leq & C \sqrt{\frac{\log(n)}{n\ell_{n}^d}} ~~~~~~~~~~
\mbox{a.s.},
\end{eqnarray*}
where $C$ is a positive constant. Therefore, we have
\begin{equation}\label{reo}
\sup_{x_{1} \in I_1} \mid \widehat  \eta_{1}(x_{1}) - \widehat
{\widehat \eta}_{1} (x_{1}) \mid  =  O\Big( \sqrt{\frac
{\log(n)}{n\ell_{n}^d}}\Big)~~~\mbox{a.s.}.
\end{equation}
Combining (\ref{reo}) and the assumption $(H.3)$, we obtain,
\begin{equation}\label{T1}
\sqrt{\frac{nh_{1,n}}{2|\log h_{1,n}|}}~~~ \sup_{x_1\in I_1}
T_{1}(x_1)=o(1) ~~~~\mbox{a.s.}
\end{equation}
Next,  turning our attention to $T_3(x_1)$. It hold that
\begin{eqnarray}
\Big|m({\bf x}) - E ( \widehat m_{n} ({\bf x}) )\Big|
 & =   & \Big| m({\bf x}) -   E \Big( \frac{Y_{i}}{ f({\bf X}_{i})} \Big[ \prod_{l=1}^{d} \frac{1}{h_{l,n}}  K_{l}\Big(\frac {x_{l} - X_{i,l}}{h_{l,n}}\Big)\Big] \Big) \Big| \nonumber \\
                                          & =   & \Big| m({\bf x}) -   \int_{\mathbb{R}^d} m({\bf u}) \Big[ \prod_{l=1}^{d} \frac{1}{h_{l,n}}   K_{l}\Big(\frac {x_{l} - u_{l}}{h_{l,n}}\Big) \Big]  d{\bf u}
                                          \Big|. \nonumber
\end{eqnarray}
The change of variables and the Taylor expansion to order $k$,
gives, with $0< \theta < 1$ and ${\bf h}=(h_{1,n},...,h_{d,n})^T$
\begin{eqnarray}
&&\hspace{-1.5cm}\Big|m({\bf x}) - E ( \widehat m_{n} ({\bf x}) )\Big| \nonumber\\& =    &\Big| m({\bf x}) - \int_{\mathbb{R}^d} m(-{\bf v h} \theta + {\bf x}) \Big( \prod_{i=1}^d    K_{i}(v_{i}) \Big) d{\bf v} \Big | \nonumber \\
                                          & \leq & \int_{\mathbb{R}^d} |(m({\bf x}) - m(-{\bf v h} \theta  + {\bf x}) ) | \Big( \prod_{i=1}^d    K_{i}(v_{i}) \Big) d{\bf v} \nonumber \\
                                          & \leq & \frac{1}{k!}\int_{\mathbb{R}^d} \sum_{k_1+...+k_d=k} h_{1,n}^{k_1} v_{1}^{k_1}...h_{d,n}^{k_d} v_{d}^{k_d}|\frac {\partial^{k} m}{\partial v_{1}^{k_1}...\partial v_{d}^{k_d}} (-{\bf v h} \theta  + {\bf x} )| \Big( \prod_{i=1}^d    K_{i}(v_{i}) \Big)  d{\bf v} \nonumber \\
                                          & \leq & \frac{1}{k!}\Big\|\frac{\partial^k m}{\partial v_1^{k_1}...\partial v_d^{k_d}}\Big\|_{\infty} \sum_{k1+...+k_d=k} h_{1,n}^{k_1}...h_{d,n}^{k_d}  \Big (\int_{\mathbb{R}} v_{1}^{k_1}...v_{d}^{k_d} \mathbb{K}({\bf v}) d{\bf v} \Big).\nonumber
\end{eqnarray}
It follows that,
\begin{eqnarray}\label{inequl}
\sup_{x_1 \in I_1} T_3(x_1) & \leq & 2 \sup_{{\bf x}\in I}\Big|
m({\bf x}) -  E(\widehat m_n ({\bf x} ))\Big| \nonumber \\
& \leq &  \frac{2}{k!}\Big\|\frac{\partial^k m}{\partial
v_1^{k_1}...\partial v_d^{k_d}}\Big\|_{\infty}
\sum_{k_1+...+k_d=k} h_{1,n}^{k_1}...h_{d,n}^{k_d}  \Big
(\int_{\mathbb{R}} v_{1}^{k_1}...v_{d}^{k_d} \mathbb{K}({\bf v})
                                          d{\bf v}\Big).
\end{eqnarray}
By combining the assumption $(H.3)$ and the statement
(\ref{inequl}), we obtain,
\begin{eqnarray}\label{T3}
\sqrt{\frac{nh_{1,n}}{2|\log h_{1,n}|}}~\sup_{x_1 \in I_1}
T_3(x_1)=o(1).
\end{eqnarray}
Finally we evaluate the term $T_2(x_1)$. From (\ref{alphag}), we
observe that,
\begin{eqnarray}
nh_1 \{\widehat {\widehat \eta}_{1}(x_{1}) - E(\widehat {\widehat
\eta}_{1}(x_{1})) \} = \sqrt{n} \alpha_{n}(g_{n}^{x_1 })
-\int_{\mathbb{R}} \sqrt{n} \alpha_{n}(g_{n}^{x_1 }) q_1 (x_1)
dx_1,
\end{eqnarray}
then
\begin{eqnarray*}
\sqrt {\frac {n h_{1,n}}{2|\log h_{1,n}|}} ~T_2(x_1) =
\frac{\pm~\alpha_{n}\Big(g_{n}^{x_1 } \Big)}{\sqrt{2h_{1,n} |\log
h_{1,n}|}} \mp \int_{\mathbb{R}} \frac{\alpha_{n}\Big(g_{n}^{x_1 }
\Big)  }{\sqrt{2h_{1,n} |\log h_{1,n}|}} q_1 (x_1 ) dx_1.
\end{eqnarray*}
We will use in this proof \cite{Camlong-Viot-Sarda-and-Vieu2000}
notations,
\begin{eqnarray*}
\hat \alpha_1 (x_1)&=& \frac{1}{nh_{1,n}}\sum_{i=1}^n\frac{\tilde
Y_{i,n}}{f_1(X_1)}K_1\Big(\frac{x_1-X_{i,1}}{h_{1,n}}\Big),\\
\tilde Y_{i,n}&= &\int_{\mathbb{R}^{d-1}}
\Big(\prod_{l=2}^d\frac{1}{h_{l,n}}K_l\Big(\frac{x_l-X_{i,l}}{h_{l,n}}\Big)\Big)\frac{q_{-1}({\bf
x}_{-1})}{f(X_{i,-1}|X_{i,1})} d{\bf x}.
\end{eqnarray*}
 Using the Cauchy-Schwartz
inequality, we obtain,
\begin{eqnarray}\label{debbarh}
& &\Big|\int_{\mathbb{R}} \frac{\alpha_{n}\Big(g_{n}^{x_1} \Big)
}{\sqrt{2h_{1,n}|\log h_{1,n}|}} ~q_1 (x_1 ) dx_1 \Big| \nonumber
\\& = & \sqrt{\frac{nh_{1,n}}{2|\log h_{1,n}|}}
\Big|\int_{\mathbb{R}^d} \{\widehat{m}_n ({\bf x}) -
E(\widehat{m}_n
({\bf x}))\} q({\bf x}) d{\bf x}\Big|  \nonumber \\
&=& \sqrt{\frac{nh_{1,n}}{2|\log h_{1,n}|}} \Big|\int_{\mathbb{R}}
\{\hat{\alpha}_1 (x_1) -
E(\hat{\alpha}_1(x_1))\} q_1( x_1) d x_1\Big| \nonumber \\
&\leq & \sqrt{\frac{nh_{1,n}}{2|\log h_{1,n}|}}
\Big[\int_{\mathcal{C}_1} \{\hat{\alpha}_1 (x_1) -
E(\hat{\alpha}_1(x_1))\}^2 dx_1\Big]^{1/2}\times\Big[
\int_{\mathbb{R}} q_1^2(x_1) dx_1\Big]^{1/2}.
\end{eqnarray}
By \cite{Camlong-Viot-Sarda-and-Vieu2000} we have, for all $x_1
\in \mathcal{C}_1$,
\begin{eqnarray*}
{\rm Var} (\hat \alpha_1(x_1))= E\Big(\hat{\alpha}_1 (x_1) -
E(\hat{\alpha}_1(x_1))
\Big)^2=\mathcal{O}\Big(n^{-2k/(2k+1)}\Big),
\end{eqnarray*}
then
\begin{eqnarray}\label{esperance}
\int_{\mathcal{C}_1}{\rm Var} (\hat \alpha_1(x_1)) dx_1 &=&
\int_{\mathcal{C}_1} E\Big(\hat{\alpha}_1 (x_1) -
E(\hat{\alpha}_1(x_1)) \Big)^2 dx_1\nonumber\\
&=& E\Big(\int_{\mathcal{C}_1} \Big(\hat{\alpha}_1 (x_1) -
E(\hat{\alpha}_1(x_1)) \Big)^2 dx_1\Big)\nonumber\\
&=& \mathcal{O}\Big(n^{-2k/(2k+1)}\Big).
\end{eqnarray}
From (\ref{esperance}), it follow that,
\begin{eqnarray}\label{sans_esp}
\int_{\mathcal{C}_1} \Big(\hat{\alpha}_1 (x_1) -
E(\hat{\alpha}_1(x_1)) \Big)^2 dx_1 =
\mathcal{O}\Big(n^{-2k/(2k+1)}\Big) ~~\mbox{a.s.}
\end{eqnarray}
By combining (\ref{debbarh}), (\ref{sans_esp}) and the assumption
$(H.4)$, we arrive at
\begin{eqnarray}\label{T2}
\sqrt {\frac {n h_{1,n}}{2|\log h_{1,n}|}} ~\sup_{x_1\in
I_1}T_2(x_1) &= &\frac{\sup_{x_1\in
I_1}\pm~\alpha_{n}\Big(g_{n}^{x_1 } \Big)}{\sqrt{2h_{1,n} |\log
h_{1,n}|}} \mp \int_{\mathbb{R}} \frac{\alpha_{n}\Big(g_{n}^{x_1 }
\Big) }{\sqrt{2h_{1,n} |\log h_{1,n}|}} q_1 (x_1 ) dx_1 \nonumber \\
& =& \frac{\sup_{x_1\in I_1}\pm~\alpha_{n}\Big(g_{n}^{x_1 }
\Big)}{\sqrt{2h_{1,n} |\log h_{1,n}|}} + o(1) ~~\mbox{a.s}.
\end{eqnarray}
Finally, note that (\ref{T1}), (\ref{T3}), (\ref{T2}) and  Lemma
\ref{lem1} are sufficient to finish the proof of Theorem
\ref{th1}.
\subsection{Proof of Theorem \ref{th2}}
Observe that,
\begin{eqnarray}
&   &\hspace{-3cm} \sqrt{\frac{nh_{1,n}}{2|\log(h_{1,n})|}}
\sup_{{\bf x} \in I} \pm \{\widehat m_{add} ({\bf x}) -
m_{add}({\bf x})
 \} - \sum_{l=1}^{d} \sigma_l , \nonumber\\
 & = &
\sum_{l=1}^{d} \Big\{ \sqrt{\frac{nh_{1,n}}{2|\log h_{1,n}|}}
\sup_{x_l\in I_l} \pm \{\widehat \eta_{l} ( x_l ) - \eta_{l}( x_l
)\}
 - \sigma_l
 \Big\}\\&& + \sqrt{\frac{nh_{1,n} }{2|\log(h_{1,n} )|}} \int_{\mathbb{R}^d}\{\widehat m_n({\bf x})- m({\bf x})\}q({\bf x})d{\bf x} \nonumber.
\end{eqnarray}
By proceeding exactly as we did along the proof of (\ref{T3}) and
(\ref{debbarh}), we arrive at, under the assumption $(H.3)$ and
$(H.4)$,
\begin{eqnarray}\label{li}
\sqrt{\frac{nh_{1,n} }{2|\log(h_{1,n} )|}}
\int_{\mathbb{R}^d}\{\widehat m_n({\bf x})- m({\bf x})\}q({\bf
x})d{\bf x}=o(1) ~~\mbox{a.s.}.
\end{eqnarray}
By combining Theorem \ref{th1} and the statement (\ref{li}), we
conclude the result of Theorem \ref{th2}.$\square$
\section{ Appendix}
Here we gather together some basic Facts, that we need for the proofs. See, for instance \cite{Mason2000} and \cite{EinmahlMason2005}.\\
\\
Let $\mathcal{G}$ be a pointwise measurable class of functions
satisfying the conditions ($\mathcal{C}$), whenever there exists a
all $ x \in \Xi$,
$$ \mathbb{G}({\bf X}) \geq \sup_{g \in \mathcal{G}} |g(x)|$$
and for some $0 < \nu, C_0 < \infty$,
$$ N(\epsilon, \mathcal{G}) < C_0 \epsilon^{-\nu}, 0 < \epsilon < 1,$$
 with
$$ N( \epsilon, \mathcal{G}) = \sup_{Q} N( \epsilon \sqrt{Q(\mathbb{G}^2)}, \mathcal{G}, d_{Q}),$$
where the supremum is taken over all probability measures $Q$ on $(\Xi, A)$ for which $ 0 < Q(\mathbb{G}^2) < \infty $ and $ d_{Q} $ is the $L_{2} ( Q) $-metric. As usual $N(\epsilon, \mathcal{G}, d)$ is the minimal number of balls $\{ g : d(g,h) < \epsilon \} $ of $d$-radius $ \epsilon $ needed to cover $\mathcal{\mathcal{G}} $.\\
\\
{\bf Fact 1.} Let $\epsilon_{1} ,....,\epsilon_{n}$ be a sequence of a random variables independent of the random vectors  $X_{1}, ....,X_{n}$. The following inequality is due to \cite{Talgrand}.\\
\\
Let $0 < M < \infty$ be a constant, such that
\begin{equation}\label{C}
|| g || \leq M, \ \ \mbox{for all} \ \ g \in \mathcal{G}.
\end{equation}
 Then for all $t
> 0$, we have  suitable finite constants $ A_{1} , A_{2} > 0,$
such that
$${\bf P} \Big\{ || n^{1/2} \alpha_{n} ||_{\mathcal{G}} \geq A_{1} \Big( E|| \sum_{i=1}^{n} \epsilon_{i} g({\bf X}_{i}) ||_{\mathcal{G}} + t \Big) \Big\} \leq 2\Big[ \exp\Big(-\frac{A_{2} t^2}{n \sigma_{\mathcal{G}}^{2}}\Big) + \exp \Big( - \frac{A_{2} t}{M}\Big)\Big],$$
where $ \sigma_{\mathcal{G}}^{2} = \sup_{g \in \mathcal{G}} Var( g({\bf X}))$.\\
\\
{\bf Fact 2.} (\cite{Mason2000}) Let $\mathcal{G}$ be a pointwise
measurable class of bounded functions such that for some constants
 $\nu,$ $ C>1$, $0 <\sigma \leq \beta $ and let $\mathbb{G}$ denote a function fulfilling the above assumptions. Assume in addition that the following four conditions hold \\
\\
${\bf A}.1$ \ $E(\mathbb{G}^2({\bf X})) \leq \beta^2.$\\
${\bf A}.2$ \ $N(\epsilon, \mathcal{G}) \leq C \epsilon^{-\nu} \ , \ 0 < \epsilon < 1$.\\
${\bf A}.3$ \ $\sup_{g \in \mathcal{G}} E[g^2({\bf X})] \leq \sigma^2.$ \\
${\bf A}.4$ \ $ \sup_{g \in \mathcal{G}} || g || \leq \frac{1}{4}\sqrt{\frac{n \sigma^2 } {\nu\log( C_1 \beta/\sigma))}}.$\\
\\
We have for a universal constant $A_{3} > 0,$\\
$$ E||\sum_{i=1}^{n} \epsilon_{i} g({\bf X}_{i}) ||_{\mathcal{G}} \leq A_{3} \sqrt{n \sigma^2\nu \log( C_1\beta / \sigma )}.$$
 \\
{\bf Fact 3.} ({\it Bernstein's inequality}) . Let $
Z_{1},....,Z_{n} $ be independent random variables with mean 0 and
identical variance $ 0  < \sigma^2 : = Var(Z_{r}) < \infty , \ r =
1,....,n $. Assume further that for some $ M > 0, |Z_{r} | < M, r=
1,....,n $. Then for all $ t > 0 $
$$ {\bf P}(Z_{1} +...+ Z_{n} > t \sqrt{n} ) \leq \exp\Big(-\frac{t^2}{2 \sigma^2 + 2/3 M n^{-1/2} t}\Big).$$
\\
{\bf Fact 4.} Let $\eta_1,...,\eta_n$ be independent mean zero
random variables with $s_n^2 = \sum_{i=1}^n E\eta_i^2 > 0$ and
$P\{|\eta_i|\leq d_i\}=1,$ where $0 < d_i \uparrow$, $1 \leq i
\leq n.$ if $\lim_{n \rightarrow \infty} d_nx_n/s_n = 0,$ where
$x_n
> x_0 > 0,$ then for every $\gamma \in (0,1),$  there is a
$C_\gamma \in (0, 1/2),$ such that for all large n
$$\mathbb{P}\Big\{\sum_{i=1}^n \eta_i \geq (1- \gamma)^2 s_n x_n \Big\} > C_\gamma \exp\Big(-x_n^2 (1- \gamma)(1- \gamma^2)/2\Big).$$
\\

\end{document}